\documentstyle[12pt]{amsart}
\newtheorem{theorem}{Theorem}[section]
\newtheorem{corollary}[theorem]{Corollary}
\newtheorem{lemma}[theorem]{Lemma}
\newtheorem{example}[theorem]{Example}

\theoremstyle{remark}
\newtheorem{remark}{Remark}

\newcommand{\Sing}{\operatorname{Sing}}
\newcommand{\Cl}{\operatorname{Cl}}
\newcommand{\orb}{\operatorname{orb}}
\newcommand{\cokernel}{\operatorname{cokernel}}
\newcommand{\Pic}{\operatorname{Pic}}
\newcommand{\Spec}{\operatorname{Spec}}
\newcommand{\closure}{\operatorname{closure}}
\newcommand{\Hom}{\operatorname{Hom}}
\newcommand{\TNemb}{\operatorname{T_Nemb}}
\newcommand{\SF}{\operatorname{SF}}
\newcommand{\Zar}{\operatorname{Zar}}
\newcommand{\ord}{\operatorname{ord}}

\begin{document}

\title{The Cohomological Brauer Group of a Toric Variety}

\author{F.R. DeMeyer}
\address{Dept. of Mathematics\\Colorado State University\\Ft. Collins, CO
80523}
\email{demeyer@@euclid.math.colostate.edu}

\author{T.J. Ford}
\address{Department of Mathematics\\Florida Atlantic University\\Boca Raton, FL
33431}
\email{Ford@@acc.fau.edu}

\author{R. Miranda}
\address{Dept. of Mathematics\\Colorado State University\\Ft. Collins, CO
80523}
\email{miranda@@riemann.math.colostate.edu}
\thanks{Second and third authors supported in part by the NSF under grants
DMS-9025092 and DMS-9104058.\\ \indent January 1992}

\maketitle

\medskip
\leftline{\bf Introduction}
\medskip

Toric varieties are a special class of rational varieties
defined by equations of the form {\it monomial = monomial}.
For a good brief survey of the history and role of toric varieties see [10].
Any toric variety $X$ contains a cover by affine open sets
described in terms of arrangements (called fans)
of convex bodies in $\Bbb R^r$.
The coordinate rings of each of these affine open sets
is a graded ring generated over the ground field by monomials.
As a consequence, toric varieties provide a good context
in which cohomology can be calculated.
The purpose of this article is to describe the second \'etale cohomology group
with coefficients in the sheaf of units of any toric variety $X$.
This is the so-called cohomological Brauer group of $X$.

Let $X = \TNemb(\Delta )$ be a toric variety which determines
and is determined by the fan $\Delta $.
To say $\Delta$ is a fan means $\Delta$ consists of finitely many elements,
each of which is a convex cone in $\Bbb R^r$,
satisfying the following conditions:
\begin{enumerate}
\item If $\sigma$, $\tau \in \Delta$, then $\sigma \cap \tau$ is a face of both
$\sigma$ and $\tau$.

\item If $\sigma \in \Delta$, then any face of $\sigma$ is in $\Delta$.

\item If $\sigma \in \Delta$ and $\vec v \in \sigma$,
then $\Bbb R\vec v \cap \sigma = {\Bbb R}_{\geq 0}\vec v \cap \sigma$
(the cones in $\Delta$ are strongly convex).

\item If $N = {\Bbb Z}^r$ is the set of lattice points in ${\Bbb R}^r$
and $\sigma \in \Delta$,
then there are $\vec v_1,\ldots ,\vec v_s \in \sigma \cap N$ with
$\sigma  = {\Bbb R}_{\geq 0}\vec v_1 + \ldots + {\Bbb R}_{\geq 0}\vec v_s$.
\end{enumerate}

We can make a fan $\Delta$ a topological space
by letting the open sets in the topology be the subsets of $\Delta$
which are themselves fans
(the open sets of $\Delta$ are the subfans of $\Delta$).
The support $|\Delta |$ of $\Delta$ is
${\bigcup_{\sigma\in\Delta} \sigma} \subseteq {\Bbb R}^r$.
A support function $h$ on $\Delta $ is a real valued function on $|\Delta |$
such that for each cone $\sigma \in \Delta $,
$h|_\sigma $ is linear on $\sigma $ and integer valued on $\sigma \cap N$.
Let $\SF(\Delta )$ be the abelian group of support functions on $\Delta $.
Associate to $\Delta $ the sheaf ${\cal S}{\cal F}$
defined by ${\cal S}{\cal F} (\Delta ') = \SF(\Delta ')$
for each open set $\Delta ' \subseteq \Delta $.
The ``restriction maps'' in the sheaf ${\cal S}{\cal F}$
are ordinary restriction.

This article compares several cohomology groups and it will be helpful to the
reader to lay out our notation with some care.
If $X$ is a irreducible scheme
and $\cal F$ is a sheaf on the \'etale site,
$H^i(X_{\acute{e}t},\cal F)$ will denote
the $i$-th \'etale cohomology group of $X$ with values in $\cal F$.
If $K$ is the function field of $X$,
then  $H^i(K/X_{\acute{e}t},\cal F)$
is the subgroup of $H^i(X_{\acute{e}t},\cal F)$
of cocycle classes split by $K$,
i.e. are trivial when restricted to the generic point of $X$.
The notation $H^i_Z(X_{\acute{e}t},\cal F)$ denotes \'etale
cohomology with supports.
Similarly, if $\cal F$ is a Zariski sheaf on $X$, then
$H^i(X_{\Zar},\cal F)$ will denote
the $i$-th Zariski cohomology group of $X$ with values in $\cal F$.
Finally, for any topological space $D$ with open cover $\cal V$,
and $\cal F$ a sheaf on $D$,
$\check{H}^i({\cal V}/X,\cal F)$ denotes the $\check{C}$ech cohomology of $X$
with respect to the cover $\cal V$ with values in $\cal F$ [11].
The direct limit over all covers is denoted by $\check{H}^i(X,\cal F)$.
The notation $H^i(X,\cal F)$ simply denotes the derived functor cohomology.
The purpose of this paper is to compare these cohomology groups,
in particular for the sheaf of units,
on a toric variety.
In this article we use the terminology and notation of [9] for toric varieties.

If $X = \TNemb(\Delta )$ is a toric variety associated to a fan $\Delta $,
$K$ is the function field of $X$,
${\cal O}^*$ is the sheaf of units of $X$ on the Zariski site
and $\Bbb G_m$ is the sheaf of units on the \'etale site,
then our main result is:

\begin{theorem}
\label{th1}
Let $X = \TNemb(\Delta )$ be the toric variety with
function field $K$ associated to the fan $\Delta $. Then
\begin{itemize}
\item[(a)] $H^1(\Delta ,{\cal S}{\cal F} ) \cong  H^2(X_{Zar},{\cal O}^*) \cong
 H^2(K/X_{\acute{e}t},\Bbb G_m)$
\item[(b)] If $\Delta '$ is a nonsingular subdivision of $\Delta $
and $\tilde X = T_Nemb(\Delta ')$, then the sequence
\[
0 \to  H^2(K/X_{\acute{e}t},\Bbb G_m) \to
H^2(X_{\acute{e}t},\Bbb G_m) \to
H^2(\tilde X_{\acute{e}t},\Bbb G_m) \to  0
\]
(with natural maps) is split exact.
\end{itemize}
\end{theorem}

By the exact sequence of Theorem \ref{th1}.b,
the calculation of the group $H^2(X_{\acute{e}t},\Bbb G_m)$ is reduced to the
calculation
of the two groups $H^2(K/X_{\acute{e}t},\Bbb G_m)$ and $H^2(\tilde
X_{\acute{e}t},\Bbb G_m)$.
The group $H^2(K/X_{\acute{e}t},\Bbb G_m)$ can be calculated in terms of
the cohomology of the finite space $\Delta $, by Theorem \ref{th1}.a.
A description of the group $H^2(\tilde X_{\acute{e}t},\Bbb G_m)$ is given in
Section 1 of [1].
In [1] Section 2, $H^2(K/X_{\acute{e}t},\Bbb G_m)$ was calculated
in terms of an integer matrix associated to the fan $\Delta $
in the case where the codimension of the singular locus
$\Sing(X)$ was no more than $2$ (see also [2]).

In order to verify our main result we first show that the cohomology of
$\Delta$ with coefficients in the sheaf ${\cal SF}$ can be computed using
$\check C$ech cohomology.
Next we show that the Zariski derived functor
cohomology for the sheaf ${\cal O}^*$ agrees with the
$\check C$ech cohomology with respect to the open cover $\{ U_{\sigma} |
\sigma \in \Delta \}$ determined by $\Delta$ (see p.6 of [9]). We then show
that for each closed point $p \in X$, the natural map $\Cl{({\cal O}_{X,p})}
\rightarrow \Cl{({\cal O}^h_{X,p})}$ from the divisor class group of the
local ring of $X$ at $p$ to the henselization is an isomorphism (see p.57
of [7]). As a consequence, $H^2(X_{\Zar},{\cal O}^*) \cong
H^2(K/X_{\acute{e}t},{\Bbb G}_m)$. The isomorphism $H^1(\Delta,{\cal SF}) \cong
H^2(X_{\Zar},{\cal O}^*)$ is explicitly defined in terms of cocycles. We
end with calculations of the groups $H^p((U_{\sigma})_{\acute{e}t},{\Bbb G}_m)$
for $U_{\sigma}$ an affine toric variety. These calculations extend those
of [1], Section 2 and as a result we show that the main result of this
article does not extend to the cohomology groups of degrees greater than 2.

\medskip
\samepage{
\leftline{\bf The cohomology of $\Delta$}
\medskip

Let $\Delta $ be a fan in ${\Bbb R}^r$.
}
Define a topology on $\Delta $ by letting the open sets in $\Delta $
be the subfans of $\Delta $.
One can think of $\Delta $ as the orbit space of $X = \TNemb(\Delta )$
under the group action of the torus $T_N$
together with the induced topology.
There is also a continuous embedding $\Delta \to X$
which sends $\sigma \in \Delta $
to the generic point of the closure of $\orb(\sigma)$.
If $\sigma \in \Delta$,
then the smallest open set containing $\sigma$ is
$\Delta (\sigma ) = \{\tau  \in  \Delta  | \tau  < \sigma \}$.
If $\{\sigma _i\}$ is the set of maximal cones in $\Delta $,
then $\{\Delta (\sigma _i)\}$ is an open cover of $\Delta $
which is a refinement of any open cover of $\Delta $ (see [4] p.223).
We call this cover the {\em finest open cover} of $\Delta $.
For any sheaf ${\cal F}$ on $\Delta $,
the $\check C$ech cohomology groups
$\check H^p(\Delta ,{\cal F} )$ can be calculated as
$\check C$ech cocycles modulo coboundaries with respect to the finest cover
$\{\Delta (\sigma _i)\}.$

\begin{lemma}
\label{lem7}
Let $\sigma $ be a cone in ${\Bbb R}^r$,
$r \geq 1$ and ${\cal F}$ any sheaf on $\Delta (\sigma )$.
Then
\begin{itemize}
\item[(a)] $H^p(\Delta (\sigma ),{\cal F} ) = 0$ for all $p \geq 1$, and
\item[(b)] if $Y$ is any non-empty closed subset of
$Z = \Delta (\sigma ) - \{0\}$,
then $H^p(Y,{\cal F} |_Y) = 0$ for all $p \geq  1.$
\item[(c)] Let $\Delta$ be a fan and ${\cal G}$ a sheaf on $\Delta $. Then
$\check H^p(\Delta ,{\cal G} ) \cong  H^p(\Delta ,{\cal G} )$
for all $p \geq 0$.  I.e., the $\check C$ech cohomology agrees with the derived
functor cohomology.
\end{itemize}
\end{lemma}

\begin{pf}
{\it Step 1.}
If $\dim \ \sigma  = 0$, then $\Delta (\sigma ) = \{0\}$
and ${\cal F}$ is a constant sheaf so (a) and (b) are both true.
Assume $\dim \ \sigma > 0$.
Let $U$ be the open set $U = \{0\} \subseteq \Delta (\sigma )$
and let $Z = \Delta (\sigma ) - U$.
We have inclusion maps $j: U \to \Delta (\sigma )$,
$i : Z \to  \Delta (\sigma )$
and an exact sequence of sheaves on $\Delta (\sigma )$:
\[
0 \to  {\cal F}_U \to  {\cal F} \to {\cal F}_Z \to  0
\]
where we write ${\cal F}_U = j_!({\cal F} |_U)$ and
${\cal F}_Z = i_*({\cal F} |_Z)$.
Since $U = \{0\}$ is the generic point of $\Delta (\sigma )$,
${\cal F}_U$ is the sheaf
\[
{\cal F}_U(V) = \left\{ \begin{array}{l}
{\cal F} (\{0\})\text{ if }V = \{0\} \\
0 \text{ if }V \neq \{0\}
\end{array}\right.
\]

Let {\it A} denote the abelian group ${\cal F} (\{0\})$
and ${\cal A}$ the constant sheaf defined by {\it A} on $\Delta (\sigma )$.
Then we have an exact sequence
\[
0 \to  {\cal F}_U \to  {\cal A} \to  {\cal A}_Z \to  0
\]
Since ${\cal A}$ is a constant sheaf
on the irreducible space $\Delta (\sigma )$, it is flasque, hence acyclic.
To prove (a), it suffices to show ${\cal A}_Z$ and ${\cal F}_Z$ are acyclic.
By [4], p.209 $H^p(\Delta (\sigma ),{\cal F}_Z) \cong  H^p(Z,{\cal F} |_Z)$,
so it suffices to prove (b).

{\it Step 2.} We prove (b) by induction on the cardinality of the set {\it Y}.
The sheaf ${\cal F} |_Y$ is the sheaf associated to the presheaf
$W \mapsto \lim_{V\cap Z=W} {\cal F}(V)$.
Given any open set $W \subseteq  Z$, there is a minimal open set
$V \subseteq  \Delta (\sigma )$ such that $V \cap  Z = W$.
Therefore the presheaf defined by $W \mapsto  {\cal F} (V)$,
where $V$ is the minimal neighborhood of $W$,
has associated sheaf ${\cal F} |_Y$.
It follows from [4] p.64 that
$\Gamma (Y,{\cal F} |_Y) = \Gamma (\Delta (\sigma ),{\cal F} )$
since the only neighborhood of $\sigma \in Y$ is the open set $Y$ itself
and the only neighborhood of $Y$ in $\Delta (\sigma )$
is $\Delta (\sigma )$ itself.
Let $\rho_1, \dots, \rho _n$ be the distinct generic points of {\it Y}.
Let $Y_i$ be the closure of $\rho _i$.
Then $Y = Y_1 \cup \dots \cup Y_n$ is the decomposition of $Y$
into irreducible components.
As a topological space,
$Y_i$ is the orbit space of the affine toric variety $\overline{orb(\rho_i)}$.
If $Y = Y_1$ consists of just one irreducible component,
we use the procedure of Step 1 to ``remove'' the generic point $\rho _1$ of $Y$
and apply our induction hypothesis to the proper closed subset
$Y - \{\rho_1\}$.
Assume now that $Y = Y_1 \cup \dots \cup  Y_n$ and $n > 1$.
Let $Z_1 = Y_1$, $Z_2 = Y_2 \cup \dots \cup Y_n$.
Then $Y = Z_1 \cup  Z_2$.
Note that $Z_1$, $Z_2$ and $Z_{12} = Z_1 \cap Z_2$
are proper closed subsets of {\it Y}.
By induction, we assume the lemma is true for $Z_1$, $Z_2$ and $Z_{12}$.
The sequence of sheaves
\[
0 \to  {\cal F}_Y \to  {\cal F}_{Z_1}\oplus {\cal F}_{Z_2}
\to  {\cal F}_{Z_{12}} \to  0
\]
on $\Delta (\sigma )$ gives rise to the long exact sequence of cohomology
\[
0 \to  H^0(Y,{\cal F} |_Y) \to  H^0(Z_1,{\cal F} |_{Z_1})
\oplus  H^0(Z_2,{\cal F} |_{Z_2}) \to  H^0(Z_{12},{\cal F} |_{Z_{12}})
\to
\]
\[
\dots \to  H^p(Y,{\cal F} |_Y) \to  H^p(Z_1,{\cal F}
|_{Z_1}) \oplus  H^p(Z_2,{\cal F} |_{Z_2}) \to  H^p(Z_{12},{\cal F}
|_{Z_{12}}) \to  \dots
\]
{}From the calculations of global sections above and this sequence,
parts (a) and (b) of the lemma follow.
Part (c) follows from part (a)
and the usual spectral sequence argument (see 4.11 p. 225 of [4]).
\end{pf}

\medskip
\samepage{
\leftline{\bf The Zariski cohomology of a toric variety}
\medskip

Let $X = U_\sigma = \TNemb(\Delta (\sigma ))$ be the affine toric variety
associated to the $s$-dimensional cone $\sigma $ in ${\Bbb R}^r$.
}
Let ${\cal F}$ be a sheaf on {\it X}.
Say that ${\cal F}$ is {\em constant on orbits} if
for any face $\tau \leq \sigma $,
${\cal F} |_{\orb(\tau)}$ is a constant sheaf on the variety $\orb(\tau)$.

\begin{lemma}
\label{lem5}
Let $X = U_\sigma $ be an affine toric variety and
${\cal F}$  a sheaf on $X$ which is constant on orbits. Then
\begin{itemize}
\item[(a)] $H^p(X_{\Zar},{\cal F} ) = 0$ for all $p \geq  1$, and
\item[(b)] if $Y$ is a closed subset of
$Z = \cup \{\orb(\tau)  | \tau \leq \sigma $, $\dim \ \tau  \geq  1 \}$
which is a union of $T_N$-orbits,
then $H^p(Y_{\Zar},{\cal F} |_Y) = 0$ for all $p \geq  1.$
\end{itemize}
\end{lemma}

\begin{pf}
If $\dim \ \sigma  = 0$, then $X$ consists of just one orbit,
the $r$-dimensional torus $T_N$, and by hypothesis
${\cal F}$  is a constant sheaf on $X$.
Since $X$ is irreducible ${\cal F}$ is flasque, hence acyclic,
and (a) is satisfied. In this case (b) is trivially true.

Assume now that $\dim \ \sigma  > 0$. Let $U = \orb(\{0\})= X-Z$.
Then ${\cal F} |_U$ is constant,
say ${\cal F} |_U = {\cal A}|_U$ where
${\cal A}$ is the constant sheaf $V \mapsto A$ on $X$
and $A$ is the abelian group ${\cal F} (U)$.
In the notation of [4], p.210, we have exact sequences of sheaves on $X$:
\begin{equation}
\label{eq7}
0 \to {\cal F}_U \to {\cal A} \to {\cal A}_Z \to  0
\end{equation}

\begin{equation}
\label{eq8}
0 \to {\cal F}_U \to {\cal F} \to {\cal F}_Z \to  0
\end{equation}

Since ${\cal A}$  is constant on the irreducible space $X$,
${\cal A}$ is acyclic.
On global sections we note that (\ref{eq7}) is
$0 \to  0 \to A \to A \to 0$ which is exact.
Since $H^p(X_{\Zar},{\cal F}_Z) = H^p(Z_{\Zar},{\cal F} |_Z)$ (p.209 [4])
we see that (a) follows from (b) applied to the sheaves
${\cal F}$ and ${\cal A}$.
We prove (b) by induction on the number of orbits
making up the closed set $Y$.
If $Y$ consists of just the orbit $\orb(\sigma) $,
the result follows from the opening paragraph of the proof.
Assuming $Y$ contains at least 2 orbits,
there is at least 1 orbit $V = \orb(\rho )$ such that
$V$ is an open subset of $Y$.
Then $Y_1 = \closure(V)$ is an irreducible closed subset
of both $Y$ and $X$.
So $({\cal F} |_Y)|_V$ is a constant sheaf.
As before, say $({\cal F} |_Y)|_V = {\cal A} |_V$
for the constant sheaf ${\cal A}$  on $Y$.
We have the exact sequence of sheaves on $Y$:
\begin{equation}
\label{eq9}
0 \to  ({\cal F} |_Y)_V \to  {\cal A} _{Y_1} \to {\cal A} _{Y_1-V} \to  0
\end{equation}

By (b) applied recursively to $Y_1$ and $Y_1-V$,
$H^p(Y_{\Zar},{\cal A} _{Y_1}) = H^p(Y_{\Zar},{\cal A} _{Y_1-V}) = 0$
for all $p \geq 1$.
Since (\ref{eq9}) is exact on global sections,
we conclude $H^p(Y_{\Zar},({\cal F} |_Y)_V) = 0$ for all $p \geq  1$.

Next consider the exact sequence
\begin{equation}
\label{eq10}
0 \to  ({\cal F} |_Y)_V \to  {\cal F} |_Y \to ({\cal F} |_Y)_{Y-V} \to  0
\end{equation}
Applying (b) recursively to $Y-V$,
$H^p(Y_{\Zar},({\cal F} |_Y)_{Y-V}) = 0$ for all $p \geq  1$
and the lemma follows from the long exact sequence associated to (\ref{eq10}).
\end{pf}

\begin{lemma}
\label{lem3}
If $X=U_\sigma$ is an affine toric variety,
then $H^2(X_{\Zar},{\cal O}^*) = 0$.
\end{lemma}

\begin{pf}
Let $K$ be the function field of $X$ and
${\cal K}^*$ the constant sheaf defined by the group of units $K^*$.
There is induced an exact sequence of sheaves
\begin{equation}
\label{eq2}
1 \to  {\cal O}^* \to  {\cal K}^* \to  {\cal C} \to  0
\end{equation}
which defines the sheaf ${\cal C}$  of Cartier divisors on $X_{\Zar}$.
There is an exact sequence
\begin{equation}
\label{eq3}
0 \to  {\cal C}  \to  {\cal W}  \to  {\cal P} \to  0
\end{equation}
where
${\cal W} = \bigoplus_{v\in X_1} {\Bbb Z}_v$
is the sheaf of Weil divisors on $X$.
By $X_1$ we mean the set of points on $X$ of codimension $1$
and for each $v \in  X_1$, $\Bbb Z_v$ is the sheaf $i_{v*}{\Bbb Z}$
where ${\Bbb Z}$ denotes the constant sheaf on $v$ defined by ${\Bbb Z}$
and $i_v : v \hookrightarrow X$.
For any open $U \subseteq  X$,
$\Gamma(U,{\cal W}) = \bigoplus_{v\in U_1} {\Bbb Z}$.

The morphism
${\cal C} \to \bigoplus_{v\in X_1} {\Bbb Z}_v$
is defined on an open $U \subseteq  X_{\Zar}$ at the presheaf level
by sending $f \in K^*$ to $\sum _{v\in U_1} ord_v(f)\cdot v$.
Because cohomology commutes with direct sums
and a constant sheaf on an irreducible space is flasque,
it follows from (2.5, p.208 [4]) that
$H^p(X_{\Zar},{\cal W} ) =
\bigoplus_{v\in X_1} H^p(X_{\Zar},{\Bbb Z}_v) =
\bigoplus_{v\in X_1} H^p(v_{\Zar},{\Bbb Z}) = 0$
and $H^p(X,{\cal K}^*) = 0$ for all $p \geq  1$.
The long exact sequence of cohomology applied to (\ref{eq2}) gives
$H^1(X_{\Zar},{\cal C} ) \cong H^2(X_{\Zar},{\cal O}^*)$
and applied to (\ref{eq3}) gives
$H^1(X_{\Zar},{\cal C} ) \cong
\cokernel\{H^0(X_{\Zar},{\cal W} ) \to  H^0(X_{\Zar},{\cal P} )\}$.
Examining the exact sequence of sheaves (\ref{eq3}) defining ${\cal P}$,
we see since ${\cal C}$ contains all locally principal divisors
that ${\cal P}$  is the sheaf associated to the presheaf
$U \to \Cl(U)/\Pic(U)$ for each open set $U \subseteq X_{\Zar}$.
The support of ${\cal P}$  is contained in the singular locus
$Z$ of $X$ (6.11, p.141 [4]).

Since $X$ is an affine toric variety,
$X = \TNemb(\Delta (\sigma ))$
where $\sigma $ is a cone in ${\Bbb R}^r$ and
$\Delta (\sigma )$ is the fan consisting of the cone $\sigma $ and its faces.
The divisor class group $\Cl(X)$ is generated by the divisors
$V(\rho ) = \overline{\orb(\rho)}$
as $\rho $ runs through the set $\Delta (\sigma )(1)$
of $1$-dimensional faces of $\sigma $.
Each of these divisors contains $\orb(\sigma)$
since each $\rho$ is a face of $\sigma$  (p.10 [8]).
If $p$ is any point in $\orb(\sigma)$,
then restriction induces maps on class groups
\[
\Cl(X) \to  \Cl({\cal O}_{X,p}) \to \Cl({\cal O}_{X,\orb(\sigma) })
\]
Each of these maps is an epimorphism. We check the composition is an
isomorphism. We have a commutative diagram
\[
\begin{array}{ccccccc}
K^* & \to & \bigoplus_{v\in X_1} {\Bbb Z}_v & \to & \Cl(X) & \to & 0 \\
\uparrow & & \uparrow & & \updownarrow  & & \\
M & \to & \bigoplus_{\rho\in \Delta(\sigma)(1)} {\Bbb Z}\cdot V(\rho) & \to &
\Cl(X) & \to & 0
\end{array}
\]
Therefore, if $D = \sum _i n_iv_i$ is a divisor on $X$
such that for each $n_i \neq  0$,
$v_i \notin  \{V(\rho ) | \rho  \in  \Delta (\sigma )(1)\}$,
then $D$ is principal.
If $U$ is an open neighborhood of the generic point of $\orb(\sigma)$,
then every divisor $V(\rho )$ has nontrivial intersection with $U$.
Any divisor $D$ on $X$ with support on $X - U$ is principal.
Therefore, $\Cl(X) = \Cl(U)$.
Since $\Cl({\cal O}_{X,\orb(\sigma)})$ is the direct limit of the $\Cl(U)$,
we have $\Cl(X) \cong  \Cl({\cal O}_{X,\orb(\sigma)})$
with the isomorphism induced by restriction.
If $\tau < \sigma $, $q \in \orb(\tau)$ and $U_\tau = \TNemb(\Delta (\tau ))$,
then $U_\tau $ is an open subset of $X$ containing $q$
and restriction induces the diagram
\[
\begin{array}{ccccc}
\Cl(X) & \stackrel{\cong}{\longrightarrow} & \Cl({\cal O}_{X,p}) &
\stackrel{\cong}{\longrightarrow} & \Cl({\cal O}_{X,\orb(\sigma) }) \\
\downarrow & & & & \downarrow \\
\Cl(U_\tau ) & \stackrel{\cong}{\longrightarrow} & \Cl({\cal O}_{X,q}) &
\stackrel{\cong}{\longrightarrow} & \Cl({\cal O}_{U_\tau ,\orb(\tau) }) \\
\downarrow & & & & \downarrow \\
0 & & & & 0
\end{array}
\]
Let $s \in  H^0(X_{\Zar},{\cal P} )$ and let $D$ be a divisor on $X$ with
$s(\orb(\sigma) ) = |D|$ in $\Cl({\cal O}_{X,\orb(\sigma) })$.
Since the sheaf maps in ${\cal P}$  are induced by restriction,
the diagram above implies $s(q) = |D|$ in $\Cl({\cal O}_{X,q})$
for any point $q \in X$.
Therefore, $H^0(X_{\Zar},{\cal P} ) = \Cl(X)$
when $X$ is an affine toric variety
and $0 = \cokernel \{H^0(X_{\Zar},{\cal W} ) \to
H^0(X_{\Zar},{\cal P} )\} = H^1(X_{\Zar},{\cal C} ) = H^2(X_{\Zar},{\cal O}^*)$
which proves the lemma.
\end{pf}

\begin{remark}
Note that in the above proof of Lemma \ref{lem3} it has been shown that
for any affine toric variety $X$,
the sheaf ${\cal P}$ is constant on orbits.
Hence $H^p(X_{\Zar},{\cal P}) = 0$ for all $p \geq 1$ by Lemma \ref{lem5}.
\end{remark}

Let $\Delta $ be a fan and
$X = \TNemb(\Delta )$ be the associated toric variety.
Let $\{\sigma _1,\dots,\sigma _n\}$ be the set of maximal cones of $\Delta$,
and ${\cal V} = \{U_{\sigma _1},\dots, U_{\sigma _n}\}$.
Then ${\cal V}$ is an open cover of $X_{\Zar}$.
We can now show that the $\check C$ech cohomology
with respect to the open cover ${\cal V}$
agrees with derived functor cohomology for the units sheaf on $X_{Zar}$.

\begin{lemma}
\label{lem6}
\begin{itemize}
\item[(a)] If $X=U_\sigma$ is an affine toric variety, then
$H^p(X_{\Zar},{\cal O}^*) = 0$ for all $p \geq  1.$
\item[(b)] If $X = \TNemb(\Delta )$ is a toric variety,
$\{\sigma _1, \dots, \sigma _n\}$ the set of maximal cones of $\Delta$,
and ${\cal V} = \{U_{\sigma _1}, \dots, U_{\sigma _n}\}$, then
$\check H^p({\cal V}\ /X,{\cal O}^*) \cong  H^p(X_{\Zar},{\cal O}^*)$
for all $p \geq 0$.
\end{itemize}
\end{lemma}

\begin{pf}
We already know $H^1(X_{\Zar},{\cal O}^*) = \Pic(X) = 0$ (p.73 of [9])
and $H^2(X_{\Zar},{\cal O}^*) = 0$ by Lemma \ref{lem3}.
{}From the proof of Lemma \ref{lem3},
$H^p(X_{\Zar},{\cal P}) \cong  H^{p+1}(X_{\Zar},{\cal C} ) \cong
H^{p+2}(X_{\Zar},{\cal O}^*)$
for $p \geq 1$.
As remarked above,
$H^p(X_{\Zar},{\cal P} ) = 0$ for all $p \geq  1$.
This proves part (a) and part (b) follows from part (a)
and the usual spectral sequence argument
(p.100 of [8]).
\end{pf}

\medskip
\samepage{
\leftline{\bf The second \'etale cohomology group of a toric variety}
\medskip

In order to verify our main result we next show (see p.57 [7])
that for each closed point $p \in X$ the natural map
$\Cl({\cal O}_{X,p}) \to \Cl({\cal O}^h_{X,p})$
from the divisor class group of the local ring of $X$ at $p$
to the henselization is an isomorphism.
}
As a consequence,
$H^2(X_{\Zar},{\cal O}^*) \cong  H^2(K/X_{\acute{e}t},\Bbb G_m)$.

\begin{lemma}
\label{lem2}
(p.57 of [7]) Let $p$ be a closed point on the toric
variety $X$. Then the natural maps
\[
\Cl({\cal O}_{X,p}) \to \Cl({\cal O}^h_{X,p}) \to
\Cl(\widehat{{\cal O}}_{X,p})
\]
on divisor class groups are isomorphisms.
\end{lemma}

\begin{pf}
There is a unique cone $\sigma $ in the fan $\Delta $ defining $X$
such that $p \in \orb(\sigma) \subseteq U_\sigma $
where $U_\sigma $ is the open set on $X$ corresponding to $\sigma $ and
$\orb(\sigma) $ is the unique closed orbit in $U_\sigma $
under the action of the torus $T_N$.
With respect to a basis choice, the coordinate ring of $U_\sigma $ is
$R = k[e(m_1), \dots, e(m_t), x^{\pm 1}_{s+1}, \dots, x^{\pm 1}_r]$
where the $e(m_i)$ are monomials in indeterminates $x_1,\dots ,x_s$
and the closed set $\orb(\sigma)$ corresponds to the ideal
$I = \langle e(m_1), \dots, e(m_t)\rangle $.
Since $p$ is a closed point contained in $\orb(\sigma)$,
the ideal in $R$ corresponding to $p$ is
$P = \langle e(m_1), \dots, e(m_t), x_{s+1}-\alpha _{s+1}, \dots,
x_r-\alpha _r\rangle $
where $\alpha _{s+1}, \dots, \alpha _r$ are non-zero elements of $k$.

Let $\rho _1, \dots, \rho _m$ be the $1$-dimensional faces of $\sigma $
and $I_1, \dots, I_m$ the height one primes in $R$
corresponding to $\orb(\rho _1), \dots, \orb(\rho_m)$.
Then $\Cl(U_\sigma ) = \Cl(R)$ is generated by the ideals
$I_1, \dots, I_m$ and each $I_j \subseteq I$.
Since each of the $I_j$ is generated by monomials we can find a monomial $y$
 such that $R(1/y) = k[x^{\pm 1}_1, \dots, x^{\pm 1}_r]$ and
$I_1,\dots, I_m$ are the only height one primes of $R$ containing $y$.
The completion $\hat R$ of $R$ with respect to the ideal $P$ is
$\hat R =
k[[e(m_1), \dots, e(m_t), x_{s+1}-\alpha _{s+1},\dots, x_r-\alpha _r]]$.
Making the substitution $y_i = x_i-\alpha _i (s+1 \leq  i \leq  r)$ gives
$\hat R(1/y) =
k[[x_1, \dots, x_s, y_{s+1}, \dots, y_r]][x^{-1}_1, \dots, x^{-1}_s]$
is factorial so $\Cl(\hat R)$ is generated by the height one primes
containing $\hat I_1, \dots, \hat I_m$.
Since $I_j$ is the ideal corresponding to $\orb(\rho_j)$
for the $1$-dimensional face $\rho _j$ of $\sigma $,
$R/I_j$ is the coordinate ring of the closure
$\overline{\orb(\rho_j)}$ of $\orb(\rho _j)$ in $U_\sigma $.
Thus, $\hat R/\hat I_j = \widehat{(R/I_j)}$
is the completion of $R/I_j$ at a closed point of $\overline{\orb(\rho_j)}$.
Since the localization of the coordinate ring of an affine toric variety
at a closed point is normal,
$\widehat{(R/I_j)}$ is a normal domain (Theorem 32, p.320 of [13]).
Thus the $\hat I_j$ are irreducible and
$\Cl(\hat R)$ is generated by $\hat I_1, \dots, \hat I_m$.
Thus the composite of natural maps
\begin{equation}
\label{eq1}
\Cl(R) \to \Cl({\cal O}_{X,p}) \to \Cl({\cal O}^h_{X,p})
\to \Cl(\widehat{{\cal O}}_{X,p})
\end{equation}
is an epimorphism.
Mori's Theorem (p.35 [3]) implies the last two maps are monomorphisms
so $\Cl({\cal O}_{X,p}) \to \Cl({\cal O}^h_{X,p})$
is an isomorphism.
\end{pf}

\begin{remark}
We showed in the proof of Lemma \ref{lem3} that the first map in (\ref{eq1})
is an isomorphism too.
\end{remark}

\begin{lemma}
\label{lem4}
If $X$ is a toric variety with function field $K$,
then
\[
H^2(X_{\Zar},{\cal O}^*) \cong  H^2(K/X_{\acute{e}t},\Bbb G_m).
\]
\end{lemma}

\begin{pf}
The proof begins as does the proof of Lemma \ref{lem3}
except we work on the \'etale site.
Let $g :\Spec\ K \hookrightarrow X_{\acute{e}t}$
be the inclusion of the generic point and
$\Bbb G_{m,K}$ the sheaf of units on $\Spec\ K$.
For any \'etale open $U \to X_{\acute{e}t}$
the group of global sections
$\Gamma (U,g_*\Bbb G_{m,K}) = R(U)^*$
where $R(U)^*$ is the group of units
in the ring of rational functions $R(U)$ on $U$.
The sheaf ${\cal C}$  of Cartier divisors is defined
by the exact sequence of sheaves:
\begin{equation}
\label{eq4}
1 \to \Bbb G_m \to g_*\Bbb G_{m,K} \to {\cal C} \to 0
\end{equation}

Consider the exact sequence of sheaves on $X_{\acute{e}t}:$
\begin{equation}
\label{eq5}
0 \to {\cal C} \to {\cal W} \to {\cal P} \to 0
\end{equation}
where ${\cal W}  = \bigoplus_{v\in X_1} \Bbb Z_v$
is the sheaf of Weil divisors.
For each $v \in X_1$, ${\Bbb Z}_v$ is the sheaf $i_{v*}{\Bbb Z}$
where ${\Bbb Z}$ denotes the constant sheaf on $v$ defined by
${\Bbb Z}$ and $i_v : v \hookrightarrow X$.
For any \'etale open $U \to X_{\acute{e}t}$,
$\Gamma (U,{\cal W} ) = \bigoplus_{v\in U_1} \Bbb Z_v$.
The morphism ${\cal C} \to \bigoplus_{v\in X_1} \Bbb Z_v$
is defined on the presheaf level by sending $f \in R(U)^*$ to
the Weil divisor $\sum _{v\in U_1}\ord_v(f)\cdot v$
and the sheaf ${\cal P}$ is defined by the sequence (\ref{eq5}).
Since $H^1(X_{\acute{e}t},{\Bbb Z}_v) = 0$ (p.106 [8]) we have
$H^1(X_{\acute{e}t},{\cal C} ) = \cokernel \{H^0(X_{\acute{e}t},{\cal W} )
\to H^0(X_{\acute{e}t},{\cal P} )\}$.
Let $Z = \Sing(X)$. On $X-Z$, ${\cal C}$  and ${\cal W}$ are isomorphic so
${\cal P}$ has its support on $Z$.
Let $p$ be a closed point of $X$.
After localizing at the geometric point $\bar p = p$, (\ref{eq5}) becomes
\begin{equation}
\label{eq6}
0 \to {\cal C}_{\bar p} \to {\cal W}_{\bar p} \to {\cal P}_{\bar p} \to  0
\end{equation}
Now ${\cal C}_{\bar p}$ is the group
of (principal) Cartier divisors of ${\cal O}^h_{X,p}$ which is
$K({\cal O}^h_{X,p})^*/({\cal O}^h_{X,p})^*$ and
${\cal W}_{\bar p}$ is the group of Weil divisors of ${\cal O}^h_{X,p}$.
Therefore ${\cal P}_{\bar p} = \Cl({\cal O}^h_{X,p})$.

Now make the same calculation as above in the Zariski topology
(see proof of Lemma \ref{lem3}).  We find
$H^1(X_{Zar},{\cal C} ) = \cokernel \{H^0(X_{Zar},{\cal W} )
\to  H^0(X_{Zar},{\cal P} )\}$.
Moreover, for each closed point $p$ on $X$,
the stalk ${\cal P}_p$ is $\Cl({\cal O}_{X,p})$.
By Lemma \ref{lem2},
the sheaves ${\cal P}$ have the same stalks
over the closed points $p$ in $X$ in both the \'etale and Zariski topologies.
Since $X$ is locally of finite type over $k$,
Remark (b) p.65 of [8] implies the stalks over all points $p$ of $X$
of the sheaves ${\cal P}$ are the same in both topologies.
Therefore these sheaves have the same groups of global sections.
Therefore $H^1(X_{\acute{e}t},{\cal C} ) \cong  H^1(X_{Zar},{\cal C} )$.
The long exact sequences of cohomology on the \'etale and Zariski sites
applied to (\ref{eq4}) and (\ref{eq2}) reduce to
\[
\begin{array}{ccccccc}
0 &  \to & H^1(X_{\acute{e}t},{\cal C} ) & \to & H^2(X_{\acute{e}t},\Bbb G_m) &
\to & H^2(\Spec(K)_{\acute{e}t},\Bbb G_{m,K}) \\
 & & \uparrow & & \uparrow & & \uparrow \\
0 & \to & H^1(X_{Zar},{\cal C} ) & \to & H^2(X_{Zar},{\cal O}^*) & \to & 0
\end{array}
\]
Therefore $H^2(X_{Zar},{\cal O}^*) \cong  H^2(K/X_{\acute{e}t},\Bbb G_m).$
\end{pf}

\medskip
\samepage{
\leftline{\bf Comparison of cohomology groups}
\medskip

Define a sheaf ${\cal S}{\cal F}$  on $\Delta $
by assigning to each open set $\Delta ' \subseteq \Delta $
the abelian group $\SF(\Delta ')$ of support functions on $\Delta '$.
}
Let $M = \Hom(N,{\Bbb Z})$ be the dual of $N$.
There is a natural map $M \rightarrow \SF(\Delta ')$
which is locally surjective.
If ${\cal M}$  denotes the constant sheaf of $M$ on $\Delta $,
then there is an exact sequence of sheaves on $\Delta$:
\begin{equation}
\label{eq11}
0 \to  {\cal U}  \to  {\cal M}  \to {\cal S}{\cal F}  \to  0
\end{equation}
where ${\cal U}$ is defined by the sequence (\ref{eq11}).
On any open $\Delta ' \subseteq  \Delta $,
${\cal U} (\Delta ') = |\Delta '|^\perp  \cap  M =
\{m \in  M | \langle m,y\rangle  = 0$ for all $y \in  |\Delta '|\}$.
Because ${\cal M}$  is flasque,
$H^p(\Delta ,{\cal M} ) = 0$ for all $p \geq  1$, so
$H^p(\Delta ,{\cal S}{\cal F} ) \cong  H^{p+1}(\Delta ,{\cal U} )$
for all $p \geq  1.$

\begin{lemma}
\label{lem8}
Let $\Delta $ be a fan and
let $X = \TNemb(\Delta )$ be the associated toric variety.
Let $\{\sigma _1, \dots , \sigma _n\}$ the set of maximal cones of $\Delta $,
and ${\cal V}  = \{U_{\sigma _1},\dots, U_{\sigma _n}\}$. Then
\[
H^p(\Delta ,{\cal U} ) \cong  \check H^p(\Delta ,{\cal U} ) \cong
\check H^p({\cal V}\ /X,{\cal O}^*) \cong  H^p(X_{\Zar},{\cal O}^*)
\]
for all $p \geq 1$.
\end{lemma}

\begin{pf}
The first isomorphism follows from Lemma \ref{lem2}.c and the third from Lemma
\ref{lem6}.b.
The constant sheaf $U \mapsto k^*$ is both acyclic and constant on orbits,
so by Lemma \ref{lem5} it follows that
$\check H^p({\cal V}\ /X,k^*) \cong  H^p(X_{\Zar},k^*) = 0$.
Therefore
$\check H^p({\cal V}\ /X,{\cal O}^*/k^*) \cong  H^p(X_{\Zar},{\cal O}^*/k^*)$
and $H^p(X_{\Zar},{\cal O}^*) \cong H^p(X_{\Zar},{\cal O}^*/k^*)$.
By comparing the $\check C$ech complexes that define the groups
$\check H^p(\Delta ,{\cal U} )$ and $\check H^p({\cal V}\ /X,{\cal O}^*/k^*)$
we see that they are equal term for term which gives the second isomorphism.
 \end{pf}

We are now in a position to prove the main theorem stated in the introduction.

\begin{pf*}{Proof of Theorem 1}
Theorem \ref{th1}(a) now follows immediately
from Lemmas \ref{lem8} and \ref{lem4}.
To prove Theorem \ref{th1}(b),
let $\Delta $ be a fan in ${\Bbb R}^r$ and
$\Delta '$ a nonsingular subdivision of $\Delta $.
Let $\tilde X = \TNemb(\Delta ') \to  X = \TNemb(\Delta )$
be the given resolution of $X$.
Since the natural map $H^2(X_{\acute{e}t},\Bbb G_m) \to
H^2(K_{\acute{e}t},\Bbb G_m)$
factors through the monomorphism $H^2(\tilde X_{\acute{e}t},\Bbb G_m) \to
H^2(K_{\acute{e}t},\Bbb G_m)$,
there is an exact sequence
\[
0 \to  H^2(K/X_{\acute{e}t},\Bbb G_m) \to H^2(X_{\acute{e}t},\Bbb G_m) \to
H^2(\tilde X_{\acute{e}t},\Bbb G_m)
\]
and it suffices to show
$H^2(X_{\acute{e}t},\Bbb G_m) \to H^2(\tilde X_{\acute{e}t},\Bbb G_m) \to 0$ is
split exact.
In Theorem 1.1 of [1] it is shown that the cup product map
$H^1(\tilde X_{\acute{e}t},{\Bbb Z}/n) \otimes H^1(\tilde X_{\acute{e}t},\mu
_n) \to H^2(\tilde X_{\acute{e}t},\mu _n)$
followed by the Kummer map
$H^2(\tilde X_{\acute{e}t},\mu _n) \to H^2(\tilde X_{\acute{e}t},\Bbb G_m)$
is surjective onto the subgroup annihilated by $n$.
Since $H^2(\tilde X_{\acute{e}t},\Bbb G_m)$ is torsion by [6], p.71,
taking the limit over all $n$ gives an exact sequence
\begin{equation}
\label{eq11a}
H^1(\tilde X_{\acute{e}t},\Bbb Q/\Bbb Z) \otimes H^1(\tilde X_{\acute{e}t},\mu)
\to H^2(\tilde X_{\acute{e}t},\Bbb G_m) \to 0.
\end{equation}

This sequence is  split exact, which follows from the description
of the cohomological Brauer group of a nonsingular toric variety
given in Theorem 1.8 of [1].
(The group  $H^2(\tilde X_{\acute{e}t} ,\Bbb G_m)$ is torsion since $\tilde X$
is nonsingular).
Indeed, the image is the subgroup of the Brauer group generated by the
symbol algebras $(m,m')_k$, and the content of Theorem 1.8 of [1]
is that every element of the torsion part of
$H^2(\tilde X_{\acute{e}t} ,\Bbb G_m)$
can be written uniquely as
$\prod_{i,j}(m_i,m_j)_{\nu_i}^{e_{ij}}$,
where the $m_i$'s are a basis for $M$, and the $\nu_i$'s divide the
invariant factors of $N/(|\Delta'|\cap N)$.
The uniqueness gives the splitting.

Consider the commutative diagram
\[
\begin{array}{ccc}
H^1(X_{\acute{e}t},{\Bbb Q/\Bbb Z}) \otimes H^1(X_{\acute{e}t},\mu) & \to &
H^2(X_{\acute{e}t},\Bbb G_m) \\
 \downarrow & & \downarrow \\
H^1(\tilde X_{\acute{e}t},{\Bbb Q/\Bbb Z}) \otimes H^1(\tilde
X_{\acute{e}t},\mu) & \to & H^2(\tilde X_{\acute{e}t},\Bbb G_m) \\
\end{array}
\]
obtained by combining (\ref {eq11a}) with its counterpart for $X$. Since
(\ref {eq11a}) splits,
it suffices to show that the natural map
$\beta:H^1(X_{\acute{e}t},\mu _n) \rightarrow  H^1(\tilde X_{\acute{e}t},\mu
_n)$
is an isomorphism for each $n$.

{}From Kummer theory, we have a commutative diagram with exact rows:
\[
\begin{array}{ccccccc}
1 \rightarrow & \Gamma (X_{\acute{e}t},\Bbb G_m) \otimes  {\Bbb Z}/n &
\rightarrow &
H^1(X_{\acute{e}t},\mu _n) & \rightarrow & _nPic\ X & \rightarrow 0 \\
 & \downarrow \alpha & & \downarrow \beta & & \downarrow \gamma & \\
1 \rightarrow & \Gamma (\tilde X_{\acute{e}t},\Bbb G_m) \otimes  {\Bbb Z}/n &
\rightarrow & H^1(\tilde X_{\acute{e}t},\mu _n) & \rightarrow & _nPic \tilde X
& \rightarrow  0
\end{array}
\]
Since $\alpha $ is clearly an isomorphism,
it suffices to show that $\gamma $ is an isomorphism.
Consider the commutative diagram
\[
\begin{array}{ccccccccc}
0 \rightarrow & {\cal U} (\Delta ) & \rightarrow & M & \rightarrow &
\SF(\Delta ) & \rightarrow & \Pic(X) & \rightarrow 0 \\
 & \downarrow  = & & \downarrow  = & & \downarrow \delta & & \downarrow
\epsilon & \\
0 \rightarrow & {\cal U} (\Delta ') & \rightarrow & M & \rightarrow &
\SF(\Delta ') & \rightarrow & \Pic(\tilde X) & \rightarrow  0
\end{array}
\]
with exact rows, which follows from sequence (\ref{eq11})
and Lemma \ref{lem8} with $p = 1$.
An element of order $n$ in $\Pic(\tilde X)$ comes from a support function
$h \in \SF(\Delta ')$ such that $n\cdot h$ is linear.
Then $h \in M \otimes {\Bbb Z}[1/n]$ hence $h$ is ${\Bbb Z}[1/n]$-valued.
But $h$ assumes ${\Bbb Z}$ values on $ |\Delta | \cap N = |\Delta '| \cap N$.
So $h \in \SF(\Delta )$ and $\gamma $ is surjective.
Since $\delta$ is injective, $\epsilon$ is injective
and $\gamma$ is an isomorphism.
\end{pf*}

\begin{corollary}
\label{cor9}
If $X$ is a nonsingular toric variety associated to a fan $\Delta $, then
$H^p(\Delta ,{\cal U} ) = H^p(X_{\Zar},{\cal O}^*) = 0$ for all $p \geq 2$.
\end{corollary}

\begin{pf}
By Lemma \ref{lem8} it suffices to show
$H^p(\Delta ,{\cal U} ) = 0$ for all $p \geq 2$.
Let ${\cal W}$  be the sheaf of $T_N$-invariant Weil divisors on $\Delta$
defined by $\Delta ' \mapsto \bigoplus_{\rho\in\Delta '(1)} \Bbb Z\cdot\rho$.
There is a monomorphism
$\SF(\Delta ') \to \bigoplus_{\rho\in\Delta '(1)} \Bbb Z\cdot\rho$
defined by
$h \mapsto \sum _{\rho \in \Delta '(1)}h(n(\rho ))\cdot \rho $
where $n(\rho )$ is a primitive element of $\rho \cap N$.
This induces the exact sequence of sheaves on $\Delta$:
\begin{equation}
\label{eq12}
0 \to {\cal S}{\cal F}  \to {\cal W} \to {\cal P} \to 0
\end{equation}
where the sheaf ${\cal P}$ is defined by (\ref{eq12}).
Since ${\cal S}{\cal F} \cong {\cal W}$  for a nonsingular fan $\Delta $,
the support of the sheaf ${\cal P}$ is contained in the set of
$\sigma \in \Delta $ which are singular.
Assuming that $\Delta $ is nonsingular, ${\cal S}{\cal F} \cong {\cal W}$.
But ${\cal W}$ is flasque, hence acyclic which proves the corollary.
  \end{pf}

\begin{remark}
If $\{\sigma_1,\dots,\sigma_s\}$ are the maximal cones of $\Delta$
then with respect to the cover $\{\Delta(\sigma_i)\}_{i=1}^s$ of $\Delta$
and $\{U_{\sigma_i}\}_{i=1}^s$ of $X$ the map
$\phi:H^1(\Delta,{\cal S}{\cal F}) \to H^2(X_{\Zar},{\cal O}^*)$
of Theorem \ref{th1}.a
is given on the $\check C$ech
cocycle level as
\[
\phi f(i,j,k) = {ef(j,k)\cdot ef(i,j) \over ef(i,k) }
\]
where $ef(i,j)$ is the monomial associated to the element $f(i,j) \in M$.
\end{remark}

We now prove some results on the higher \'etale cohomology groups which extend
some results found in the beginning of Section 2 of [1].
\begin{corollary}
\label{lem10}
Let $\sigma $ be an $s$-dimensional cone in ${\Bbb R}^r$
and $X = \TNemb(\Delta (\sigma )) = U_\sigma $
the affine toric variety associated to $\sigma $. Then
\begin{itemize}
\item[(a)] For all $p \geq  1$, $H^p(X_{\acute{e}t},\Bbb G_m)$ is torsion.
\item[(b)] $H^p(X_{\acute{e}t},\Bbb G_m) \cong  H^p(T^{r-s}_{\acute{e}t},\Bbb
G_m)$ for all $p \geq 0$,
where $T^{r-s}$ is the $(r-s)$-dimensional torus
$\Spec k[\sigma ^\perp \cap M].$
\item[(c)] $H^p(T_N/X,G_m) = \ker( H^p(X_{\acute{e}t},\Bbb G_m) \to
H^p({(T_N)}_{\acute{e}t},\Bbb G_m) ) = 0$ for all $p \geq  0.$
\end{itemize}
\end{corollary}

\begin{pf}
Let $\bar \sigma $ be the cone $\sigma $ viewed as a cone in the
$s$-dimensional ${\Bbb R}$-vector space ${\Bbb R}\sigma $.
Then $X = U_{\bar \sigma}\times T^{r-s}$.
We can give the coordinate ring $k[{\cal S}_\sigma ]$ of
$X$ a ${\Bbb Z}_{\geq 0}$-grading such that the degree-$0$ subring
is the coordinate ring of the $T^{r-s}$-factor.
It is now clear that (b) is true for $p = 0$.
For $p = 1$, it is known that the Picard group of $X$ is trivial [9].
To finish the proof of (b), by [5] it suffices to show that
$H^p(X_{\acute{e}t},\Bbb G_m)$ is isomorphic to $H^p({(X\times \Bbb
A^1)}_{\acute{e}t},\Bbb G_m)$ by the natural map.
By the argument of [5] p.164 it suffices to prove (a).

Let $Z = X - T_N$ denote the complement of the torus in $X$.
Denote by $(X^h,Z^h)$ the henselization of the couple $(X,Z)$.
By excision (p.92 of [8]),
we have
\[
H^p_Z(X_{\acute{e}t},\Bbb G_m) \cong  H^p_{Z^h}(X^h_{\acute{e}t},\Bbb G_m)
\]
for all $p \geq 1$.
By Lemmas \ref{lem3} and \ref{lem4}
$H^2(X_{\acute{e}t},\Bbb G_m) \to H^2({(X-Z)}_{\acute{e}t},\Bbb G_m)$ is
injective.
By [6], p.71,
$H^p({(X-Z)}_{\acute{e}t},\Bbb G_m)$ and $H^p({(X^h-Z^h)}_{\acute{e}t},\Bbb
G_m)$ are torsion for all $p \geq 2$
since $X-Z$ and $X^h-Z^h$ are nonsingular.
By [12] $H^p(X^h_{\acute{e}t},\Bbb G_m) \cong H^p(Z^h_{\acute{e}t},\Bbb G_m)
\cong H^p(Z_{\acute{e}t},\Bbb G_m)$ for all $p \geq 1$.
By the commutative diagram

{\small
\[
\begin{array}{ccccccccc}
\to & H^{p-1}({(X-Z)}_{\acute{e}t},\Bbb G_m) & \to & H^p_Z(X_{\acute{e}t},\Bbb
G_m) & \to & H^p(X_{\acute{e}t},\Bbb G_m) & \to & H^p({(X-Z)}_{\acute{e}t},\Bbb
G_m) & \to \\
 & \downarrow & & \downarrow\cong & & \downarrow & & \downarrow & \\
\to & H^{p-1}({(X^h-Z^h)}_{\acute{e}t},\Bbb G_m) & \to &
H^p_{Z^h}(X^h_{\acute{e}t},\Bbb G_m) & \to & H^p(X^h_{\acute{e}t},\Bbb G_m) &
\to & H^p({(X^h-Z^h)}_{\acute{e}t},\Bbb G_m) & \to \\
\end{array}
\]
}
it suffices to show $H^p(Z_{\acute{e}t},\Bbb G_m)$ is torsion for $p \geq 3$.
Suppose $1 \leq t \leq s$ and let $\rho _1,\dots, \rho _n$
be the faces of $\sigma $ of dimension $t$.
Let $U_i = \orb(\rho _i)$, $U = U_1 \cup \dots \cup U_n$,
$Y_i = \closure(U_i)$  and $Y = Y_1 \cup \dots \cup Y_n$.
Then $Y = \cup \{ \orb(\tau) | \tau \leq \sigma, \dim\,\tau \geq  t \}$
is a closed subset of $Z$, and if $t = 1$, $Y$ is equal to $Z$.
If $t = s$, $Y = \orb(\sigma)$.
Notice that
$Y - U = \cup \{ \orb(\tau) | \tau \leq \sigma, \dim \, \tau \geq  t+1 \}$.
We use descending induction on $t$ to show that
$H^p(Y_{\acute{e}t},\Bbb G_m)$ is torsion for all $p \geq 3$.
If $t = s$, then $Y = \orb(\sigma)  \cong  T^{r-s}$.
It follows from [6], p.71 that $H^p(T^{r-s}_{\acute{e}t},\Bbb G_m)$ is torsion
for all $p \geq  2$.
Suppose $1 \leq  t < s$ and that the result is true for $t+1$.
In the above notation, we are assuming for all $p \geq  3$ that
$H^p({(Y-U)}_{\acute{e}t},\Bbb G_m)$ is torsion and we are going to prove that
$H^p(Y_{\acute{e}t},\Bbb G_m)$ is torsion.
Write $Z_1 = Y-U$.
Let $(Y^h,Z^h_1)$ denote the henselization of the couple $(Y,Z_1)$.
Since $Y-Z_1 = U$ and $Y^h-Z^h_1$ are nonsingular
$H^p(U_{\acute{e}t},\Bbb G_m)$ and $H^p({(Y^h-Z^h_1)}_{\acute{e}t},\Bbb G_m)$
are torsion for all $p \geq 2$.
By [12] $H^p(Y^h_{\acute{e}t},\Bbb G_m) \cong  H^p({(Z^h_1)}_{\acute{e}t},\Bbb
G_m) = H^p({(Z_1)}_{\acute{e}t},\Bbb G_m)$
which is torsion by induction, for all $p \geq  3$.
By excision $H^p_{Z_1}(Y^h_{\acute{e}t},\Bbb G_m) \cong
H^p_{Z^h_1}(Y^h_{\acute{e}t},\Bbb G_m)$.
The commutative diagram
{\small
\[
\begin{array}{ccccccccc}
\to & H^{p-1}(U_{\acute{e}t},\Bbb G_m) & \rightarrow &
H^p_{Z_1}(Y_{\acute{e}t},\Bbb G_m) &
\rightarrow & H^p(Y_{\acute{e}t},\Bbb G_m) & \rightarrow &
H^p(U_{\acute{e}t},\Bbb G_m) & \to  \\
  & \downarrow & & \downarrow \cong & & \downarrow & & \downarrow & \\
 \to & H^{p-1}({(Y^h-Z^h_1)}_{\acute{e}t},\Bbb G_m) & \rightarrow &
H^p_{Z^h_1}(Y^h_{\acute{e}t},\Bbb G_m) & \rightarrow &
H^p(Y^h_{\acute{e}t},\Bbb G_m) & \rightarrow &
H^p({(Y^h-Z^h_1)}_{\acute{e}t},\Bbb G_m) & \to \end{array}
\]
}
shows that $H^p(Y_{\acute{e}t},\Bbb G_m)$ is torsion for each $p \geq  3$
which proves (a) and (b).
To prove (c), note that the closed immersion $T^{r-s} \rightarrow X$
restricts to a closed immersion $T^{r-s} \rightarrow T_N$ with a section.
Therefore the natural map $H^p(X_{\acute{e}t},\Bbb G_m) \rightarrow
H^p({(T_N)}_{\acute{e}t},\Bbb G_m)$
factors into the composition of monomorphisms:
$H^p(X_{\acute{e}t},\Bbb G_m) \rightarrow  H^p(T^{r-s}_{\acute{e}t},\Bbb G_m)
\rightarrow  H^p({(T_N)}_{\acute{e}t},\Bbb G_m)$.
\end{pf}

As was shown in Theorem \ref{th1} above,
$H^2(X_{\Zar},{\cal O}^*) \cong  H^2(K/X_{\acute{e}t},\Bbb G_m)$
for any toric variety $X$ with function field $K$.
There is always an isomorphism
$H^1(X_{\Zar},{\cal O}^*) \cong  H^1(K/X_{\acute{e}t},\Bbb G_m) \cong \Pic(X)$,
but the next example shows that for $p > 2$ the natural map
$H^p(X_{\Zar},{\cal O}^*) \rightarrow H^p(K/X_{\acute{e}t},\Bbb G_m)$ will not
be surjective in general.

\samepage{
\begin{example}
\label{ex13}
\mbox{ }
\end{example}
Let $\rho _0 = (1,1,1)$, $\rho _1 = (1,0,0)$, $\rho _2 = (0,1,0)$,
}
$\rho _3 = (0,0,1)$,
$\sigma _1 = {\Bbb R}_{\geq 0}\rho _0 +
{\Bbb R}_{\geq 0}\rho _1 + {\Bbb R}_{\geq 0}\rho _2$,
$\sigma _2 = {\Bbb R}_{\geq 0}\rho _0 + {\Bbb R}_{\geq 0}\rho _2 + {\Bbb
R}_{\geq 0}\rho _3$,
$\sigma _3 = {\Bbb R}_{\geq 0}\rho _0 + {\Bbb R}_{\geq 0}\rho _1 +
{\Bbb R}_{\geq 0}\rho _3$,
$\Delta  = \Delta (\sigma _1) \cup \Delta (\sigma _2) \cup \Delta (\sigma _3)$
and $X = \TNemb(\Delta )$.
Let ${\cal V}  = \{ U_{\sigma _1}, U_{\sigma _2}, U_{\sigma _3}\}$.
{}From Corollary \ref{lem10}(b) above, $H^p({(U_{\sigma _i})}_{\acute{e}t},\Bbb
G_m) = 0$ for $p > 0$, therefore $\check H^p({\cal V}\ /X_{\acute{e}t},\Bbb
G_m) = 0$ for $p > 0$.
Using either the spectral sequence
$\check H^p({\cal V}\ /X_{\acute{e}t},H^q(\Bbb G_m)) \Rightarrow
H^{p+q}(X_{\acute{e}t},\Bbb G_m)$
or a series of Mayer-Vietoris arguments ([8], p.110) one can compute:

$H^p(X_{\acute{e}t},\Bbb G_m) = \left\{\begin{array}{cc}
k^* & i=0 \\
\Bbb Z & i=1 \\
0 & i=2,3 \\
\Bbb Q/\Bbb Z & i=4 \\
0 & i \geq 5
\end{array}\right.$

{}From [6], p.71, $H^p((T_N)_{\acute{e}t},{\Bbb G}_m)$ is torsion for $p \ge
2$.
{}From [8], p.253, $H^p((T_N)_{\acute{e}t},\mu_n) = 0$ for $p \ge 4$ and $n >
1$,
where $\mu_n$ denotes the sheaf of $n$-th roots of unity. From Kummer
theory [8], p.125, $H^p((T_N)_{\acute{e}t},{\Bbb G}_m) = 0$ for $p \ge 4$.
Since
the natural map $H^p(X_{\acute{e}t},{\Bbb G}_m) \rightarrow
H^p(K_{\acute{e}t},{\Bbb
G}_m)$ factors through $H^p((T_N)_{\acute{e}t},{\Bbb G}_m)$, we have
$H^p(K/X_{\acute{e}t},{\Bbb G}_m) = H^p(X_{\acute{e}t},{\Bbb G}_m)$ is non-zero
for $p
= 4$. Therefore
\[
H^4(K/X_{\acute{e}t},{\Bbb G}_m) \neq
\check{H}^4({\cal V}/X_{\acute{e}t},{\Bbb G}_m) =
\check{H}^4({\cal V}/X_{\Zar},{\cal O}^*) =
H^4(X_{\Zar},{\cal O}^*).
\]
So Theorem 1 does not extend to $H^4$.

\medskip
\centerline{REFERENCES}

[1] F. DeMeyer and T. Ford, {\it On the Brauer group of toric varieties}, to
appear in:   Trans. Amer. Math. Soc.

[2] F. DeMeyer and T. Ford, {\it Nontrivial, locally trivial Azumaya algebras},
in: ``Azumaya Algebras, Actions, and Modules'', D. Haile and J. Osterburg,
eds., Contemporary Math. Vol 124,  Amer. Math. Soc. (1992), 39-49.

[3] R. Fossum, ``The Divisor Class Group of a Krull Domain'', Springer Verlag,
New York, Berlin, 1973.

[4] R. Hartshorne, ``Algebraic Geometry'', Springer Verlag, New York, Berlin,
1977.

[5] R. Hoobler, {\it Functors of graded rings}, in: ``Methods in Ring Theory'',
F. van Oystaeyen, ed., NATO ASI Series (Reidel, Dordrecht, 1984) 161-170.

[6] A. Grothendieck, {\it Le groupe de Brauer II}, in: ``Dix Expos\'es sur la
Cohomologie   des Sch\'emas'', North Holland, Amsterdam, 1968.

[7] G. Kempf, F. Knudsen, D. Mumford, B. Saint-Donat, ``Toroidal Embeddings
I'',
Lecture Notes in Mathematics Vol. 339, Springer-Verlag, New York, Berlin, 1973.

[8] J. Milne, ``Etale Cohomology'', Princeton University Press, Princeton,
N.J.,
1980.

[9] T. Oda, ``Convex Bodies and Algebraic Geometry'', Springer-Verlag,
Berlin, Heidelberg, 1988.

[10] M. Reid, review of ``Convex Bodies and Algebraic Geometry'' by T. Oda,
Bull. Amer. Math. Soc. {\bf 21} (1989) 360-364.

[11] J. P. Serre, {\it Faisceaux alg}\'e{\it briques coh}\'e{\it rents}, Ann.
of Math. {\bf 61} (1955) 197-278.

[12] R. Strano, {\it On the} \'e{\it tale cohomology of Hensel rings}, Comm.
Algebra, {\bf 12} (1984)   2195-2211.

[13] O. Zariski and P. Samuel, ``Commutative Algebra, Vol. II'',
Springer-Verlag,   New York, Heidelberg, 1960.

\end{document}